\documentclass[12pt]{article}
\usepackage{amssymb,amsmath,amsthm,url,textcomp,color}

\theoremstyle{plain}

\theoremstyle{definition}

\newcommand\blfootnote[1]{%
  \begingroup
  \renewcommand\thefootnote{}\footnote{#1}%
  \addtocounter{footnote}{-1}%
  \endgroup
}

\newcommand{\bF}{\mathbb{F}}

\newcommand{\GL}{\operatorname{GL}}
\newcommand{\Alg}{\operatorname{Alg}}

\newcommand{\hG}{\widehat{S}_n}
\newcommand{\hg}{\widehat{g}}
\newcommand{\hh}{\widehat{h}}

\newcommand{\ta}{\tilde{a}}
\newcommand{\tc}{\tilde{c}}

\renewcommand{\phi}{\varphi}

\begin{document}
\title{A Practical Cryptanalysis of the Algebraic Eraser}
\author{Adi Ben-Zvi\thanks{
Department of Mathematics,
Bar-Ilan University,
Ramat Gan 5290002,
Israel}
\and
Simon R.\ Blackburn\thanks{
Department of Mathematics,
Royal Holloway University of London,
Egham, Surrey TW20 0EX,
United Kingdom}
\and
Boaz Tsaban\footnotemark[1]
}
\maketitle

\begin{abstract}
We present a novel cryptanalysis of the Algebraic Eraser primitive. This key agreement scheme, based on techniques from permutation groups, matrix groups and braid groups, is proposed as an underlying technology for ISO/IEC 29167-20, which is intended for authentication of RFID tags. SecureRF, the company owning the trademark Algebraic Eraser, markets it as suitable in general for lightweight environments such as RFID tags and other IoT applications. Our attack is practical on standard hardware: for parameter sizes corresponding to claimed 128-bit security, our implementation recovers the shared key using less than 8 CPU hours, and less than 64MB of memory.
\end{abstract}

\blfootnote{\copyright\  IACR 2016. This article is the final version submitted by the author(s) to the IACR and to Springer-Verlag on 1 June 2016. The version published by Springer-Verlag is available at [DOI to follow].}

\section{Introduction}
\label{sec:introduction}

The Algebraic Eraser\texttrademark\  is a key agreement scheme using techniques from non-commutative group theory.
It was announced by Anshel, Anshel, Goldfeld and Lemieaux in 2004; the corresponding paper~\cite{AnshelAnshelGoldfeldLemieux06}  appeared in 2006. 
The Algebraic Eraser is defined in a very general fashion: various algebraic structures (monoids, groups and actions) need to be specified in order to be suitable for implementation. Anshel \emph{et al.}\ provide most of this extra information, and name this concrete realisation of the Algebraic Eraser the \emph{Colored Burau Key Agreement Protocol (CBKAP)}. This concrete representation involves a novel blend of finite matrix groups and permutation groups with infinite braid groups. A company, SecureRF, owns the trademark to the Algebraic Eraser, and is marketing this primitive as suitable for low resource environments such as RFID tags and Internet of Things (IoT) applications. The primitive is proposed as an underlying technology for ISO/IEC~29167-20, and work on this standard is taking place in ISO/IEC~JTC~1/SC~31/WG 7. The company has also presented the primitive to the Internet Research Task Force's Crypto Forum Research Group (IRTF CFRG), with a view towards standardisation. IoT is a growth area, where current widely-accepted public key techniques struggle to operate due to tight efficiency constraints. It is likely that solutions which are efficient enough for these applications will become widely deployed, and the nature of these applications make system changes after deployment difficult. Thus, it is vital to scrutinise the security of primitives such as the Algebraic Eraser early in the standardisation process, to ensure only secure primitives underpin standardised protocols.

In a presentation to the NIST Workshop in Lightweight Cryptography in 2015, SecureRF claims a security level of $2^{128}$ for their preferred parameter sizes, and compares the speed of their system favourably with an implementation of a key agreement protocol based on the NIST recommended~\cite{NIST13} elliptic curve K-283. The company reports~\cite{AtkinsGunnells15} a speed-up by a factor of 45--150, compared to elliptic curve key agreement at 128-bit security levels. It claims that the computational requirements of the Algebraic Eraser scales linearly with the security parameter, in contrast to the quadratic scaling of elliptic-curve-based key agreement.

\paragraph{Related work} 
The criteria for choosing some global parameters of the scheme (namely certain subgroups $C$ and $D$ of matrices over a finite field, and certain subgroups $A$ and $B$ of a certain infinite semidirect product of groups) are not given in~\cite{AnshelAnshelGoldfeldLemieux06}, and have not been made available by SecureRF. In the absence of this information, it is reasonable to proceed initially with a cryptanalysis on the basis that these parameters are chosen in a generic fashion. All previous cryptanalyses have taken this approach.

Myasnikov and Ushakov~\cite{MyasnikovUshakov09} provide a heuristic length-based attack on the CBKAP that works for the parameter sizes originally suggested~\cite{AnshelAnshelGoldfeldLemieux06}. However, Gunnells~\cite{Gunnells11} reports that this attack quickly becomes unsuccessful as parameter sizes grow; the parameter sizes proposed by SecureRF make this attack impractical.\footnote{There is an analogy with the development of RSA here: the size of primes (200 digits) proposed in the original article~\cite{RSA} was made obsolete by improvements in integer factorisation algorithms~\cite{BahrBoehm05}.}
Kalka, Teicher and Tsaban~\cite{KalkaTeicherTsaban12} provide an efficient cryptanalysis of the CBKAP for 
arbitrary parameter sizes. The attack uses the public key material of Alice and the messages exchanged between Alice and Bob to derive an equivalent to the secret random information generated by Bob, which then compromises the shared key, and so renders the scheme insecure. In particular, the techniques of~\cite{KalkaTeicherTsaban12} will succeed when the global parameters are chosen generically.

SecureRF uses proprietary distributions for global parameters, so the cryptanalysis of~\cite{KalkaTeicherTsaban12} attack does not imply that the CBKAP as implemented is insecure.\footnote{The analogy with RSA continues: factorisation of a randomly chosen integer $n$ is much easier than when $n$ is a product of two primes of equal size, which is why the latter is used in RSA.} Indeed, Goldfeld and Gunnells~\cite{GoldfeldGunnells12} show that by choosing the subgroup $C$ carefully one step of the attack of~\cite{KalkaTeicherTsaban12} does not recover the information required to proceed, and so this attack does not succeed when parameters are generated in this manner.

\paragraph{Our contribution}
There are no previously known attacks on the CBKAP for the proposed parameter sizes, provided the parameters are chosen to resist the attack of~\cite{KalkaTeicherTsaban12}. The present paper describes a new attack on the CBKAP that does not assume any structure on the subgroup~$C$. Thus, a careful choice of the subgroup $C$ will have no effect on the applicability of our attack, and so the proposed security measure offered by Goldfeld and Gunnells~\cite{GoldfeldGunnells12} to the attack of~\cite{KalkaTeicherTsaban12} is bypassed.

The earlier cryptanalyses of CBKAP~(\cite{MyasnikovUshakov09},\cite{KalkaTeicherTsaban12})
attempt to recover parts of Alice's or Bob's secret information.
The attack presented here recovers the shared key directly from Alice's public key 
and the messages exchanged between Alice and Bob.

SecureRF have kindly provided us with sets of challenge parameters of the full recommended size, and our implementation is successful in recovering the shared key in all cases. Our (non-optimised) implementation recovers the common key in under 8~hours of computation, and thus the security of the system is \emph{much} less than the $2^{128}$ level claimed for these parameter sizes. The attack scales well with size, so increasing parameter sizes will not provide a solution to the security problem for the CBKAP.

\paragraph{Conclusion and recommendation}
Because our attack efficiently recovers the shared key of the CBKAP for recommended parameter sizes, using parameters provided by SecureRF, we believe the results presented here cast serious doubt on the suitability of the Algebraic Eraser for the applications proposed. We recommend that the primitive in its current form should not be used in practice, and that full details of any revised version of the primitive should be made available for public scrutiny in order to ensure a rigorous security analysis.

\paragraph{Recent developments}
Since the first version of this paper was posted, there have been two recent developments. Firstly, authors from SecureRF have 
pos\-ted~\cite{benzvi_rebuttal_arxiv} a response to our attack, concentrating in the main on the implications for the related ISO standard and providing some preliminary thoughts on how they might redesign the primitive. Until the details are finalised, it is too soon to draw any conclusions on the security of any redesigned scheme, though there have already been some discussions on Cryptography Stack Exchange~\cite{stackexchange}.   Secondly, Blackburn and Robshaw~\cite{BlackburnRobshaw} have posted a paper that cryptanalyses the ISO standard itself, rather then the more general underlying Algebraic Eraser primitive. 

\paragraph{Structure of the paper}
The remainder of the paper is organised as follows. Sections~\ref{sec:notation} and~\ref{sec:protocol} establish notation, and describe the CBKAP. We describe a slightly more general protocol than the CBKAP, as our attack naturally generalises to a larger setting. We describe our attack in Section~\ref{sec:attack}. In Section~\ref{sec:implementation} we describe the results of our implementations and provide a short conclusion.

\section{Notation}
\label{sec:notation}

This section establishes notation for the remainder of the paper. We closely follow the notation from~\cite{KalkaTeicherTsaban12}, which is in turn mainly derived from the notation in~\cite{AnshelAnshelGoldfeldLemieux06}, though we do introduce some new terms.

Let $\bF$ be a finite field of small order (e.g., $|\bF|=256$) and let $n$ be a positive integer (e.g., $n=16$). 
Let $S_n$ be the symmetric group on the set $\{1,2,\ldots,n\}$, 
and let $\GL_n(\bF)$ be the group of invertible $n\times n$ matrices 
with entries in $\bF$.

Let $M$ be a subgroup of $\GL_n(\bF(t_1,\ldots,t_n))$, where the elements $t_i$ are algebraically independent commuting indeterminates. Indeed, we assume that the group $M$ is contained in the subgroup of 
$\GL_n(\bF(t_1,\ldots,t_n))$ of matrices whose determinant can be written as 
$a\mathbf{t}$ for some non-zero element $a\in\bF$ and some, possibly empty, word 
$\mathbf{t}$ in the elements $t_i$ and their inverses. Let $\overline{M}$ be the subgroup of $\GL_n(\bF(t_1,\ldots,t_n))$ generated by permuting the indeterminates of elements of $M$ in all possible ways.

Fix non-zero elements $\tau_1,\ldots,\tau_n\in \bF$. Define the homomorphism $\phi\colon \overline{M}\rightarrow \GL_n(\bF)$ 
to be the evaluation map, computed by replacing each indeterminate $t_i$ by the corresponding element $\tau_i$. Our assumption on the group $M$ means that $\phi$ is well defined.

The group $S_n$ acts on $\overline{M}$ by permuting the indeterminates $t_i$. Let $\overline{M}\rtimes S_n$ be the semidirect product of $\overline{M}$ and $S_n$ induced by this action. More concretely, if we write ${}^ga$ for the action of an element $g\in S_n$ on an element $a\in \overline{M}$, then the elements of $\overline{M}\rtimes S_n$ are pairs $(a,g)$ with $a\in \overline{M}$ and $g\in S_n$, and group multiplication is given by
\[
(a,g)(b,h)=(a\,{}^gb,gh)
\]
for all $(a,g),(b,h)\in \overline{M}\rtimes S_n$.

Let $C$ and $D$ be subgroups of $\GL_n(\bF)$ that \emph{commute elementwise}: 
$cd=dc$ for all $c\in C$ and $d\in D$. The CBKAP specifies that $C$ is a subgroup consisting of all invertible matrices of the form $\ell_0+\ell_1\kappa+\cdots +\ell_r\kappa^r$ where $\kappa$ is a fixed matrix, $\ell_i\in\bF$ and $r\geq 0$. So $C$ is the group of units in the $\bF$-algebra generated by $\kappa$. Moreover, the CBKAP specifies that $D=C$. But we do not assume anything about the forms of $C$ and $D$ in this paper, other than the fact that they commute.

Let $\Omega=\GL_n(\bF)\times S_n$ and let $\hG=\overline{M}\rtimes S_n$.
We have two actions on $\Omega$. 
Firstly, there is the right action of the group $\hG$ on $\Omega$ via a map 
\[
*\colon\Omega\times \hG\rightarrow \Omega,
\]
as defined in~\cite{AnshelAnshelGoldfeldLemieux06,KalkaTeicherTsaban12}. So
\[
(s,g)*(b,h)=(s\phi(^gb),gh)
\]
for all $(s,g)\in \Omega$ and all $(b,h)\in \hG$.
Secondly, there is a left action of the group $\GL_n(\bF)$ on $\Omega$ via the map
\[
\bullet\colon \GL_n(\bF)\times \Omega\rightarrow \Omega
\]
given by matrix multiplication:
\[
x\bullet (s,g)=(xs,g)
\]
for all $x\in \GL_n(\bF)$ and all $(s,g)\in \Omega$. Note that for all 
$x\in \GL_n(\bF)$, all $\omega\in \Omega$ and all $\hg\in \hG$ we have that
\[
(x\bullet \omega)*\hg=x\bullet(\omega*\hg).
\]
Also note that the left action is $\bF$-linear, in the sense that if $x\in \GL_n(\bF)$ 
can be written in the form
\[
x=\sum_{i=1}^r\ell_ic_i
\]
for some $c_i\in \GL_n(\bF)$ and $\ell_i\in \bF$, then for all $(s,g)\in\Omega$ we have
\[
x\bullet(s,g)=\sum_{i=1}^r\ell_i(c_i\bullet (s,g)).
\]
To interpret the right hand side of the equality above: the subset of $\Omega$ whose second component is a fixed element of $S_n$ is naturally an $\bF$-vector space, where addition and scalar multiplication takes place in the first component only.

Finally, let $A$ and $B$ be subgroups of $\hG$ that \emph{$*$-commute}: for all $(a,g)\in A$, $(b,h)\in B$ and $\omega\in\Omega$,
\[
(\omega*(a,g))*(b,h)=(\omega*(b,h))*(a,g).
\]

\section{The CBKAP protocol}
\label{sec:protocol}

\subsection{Overview}

The CBKAP is unusual in that the parties executing it, Alice and Bob,
use different parts of the public key in their computations: neither party needs to know all of the public key. The security model assumes that one party's public key material is known to the adversary: say Alice's public key material is known, but Bob's `public' key (which is better thought of as part of his private key material) is not revealed. The adversary, Eve, receives just Alice's public information, and the messages sent over the insecure channel. Security means that Eve cannot feasibly compute any significant information about~$K$. The attack in~\cite{KalkaTeicherTsaban12} works in this model. The same is true for the attack we describe below.

In a typical proposed application, the protocol might be used to enable a low-power device, such as an RFID tag, to communicate with a central server.  Data on an RFID tag is inherently insecure, as is system-wide data. So the above security model is realistic (and conservative) for these application settings. 

\subsection{The protocol}

Public parameters (for Alice) include the parameters $n$, $\bF$, $M$, $\tau_1,\ldots,\tau_n$, $C$ and $A$. The groups $M$, $C$ and $A$ are specified by their generating sets. For efficiency reasons, the generators of the group $A$ are written as words in a certain standard generating set for the group $\hG$.
We discuss this further in Section~\ref{sec:implementation}, but see the TTP algorithm in~\cite{AnshelAnshelGoldfeldLemieux06} for full information. It is assumed that Eve knows the parameters $n$, $\bF$, $M$, $\tau_1,\ldots,\tau_n$, $C$ and $A$.  Bob needs to know the groups $B$ and $D$, rather than the groups $A$ and $C$. Eve does not need to know the subgroups $B$ and $D$ for our attack to work.

We write $e$ for the identity element of $S_n$. We write $I_n$ for the identity matrix in 
$\GL_n(\bF)$, and write $1=(I_n,e)\in \Omega$.

Alice chooses elements $c\in C$ and $\hg=(a,g)\in A$. She computes the product
\[
c\bullet 1 * \hg=(c\phi(a),g)\in \Omega
\]
and sends it to Bob over an insecure channel.

Bob, who knows the groups $B$ and $D$, chooses elements $d\in D$ and $\hh=(b,h)\in B$. He computes the product
\[
d\bullet 1 *\hh=(d\phi(b),h)\in \Omega
\]
and sends it to Alice over the insecure channel.

Note that $cd=dc$ because $c\in C$ and $d\in D$, and the groups $C$ and $D$ commute
elmentwise. Thus,
\begin{align*}
d\bullet (c\bullet 1 * \hg)*\hh&=(dc)\bullet (1 *\hg)*\hh\\
&=(cd)  \bullet (1 *\hg)*\hh\\
&=(cd)\bullet (1*\hh)*\hg\\
&\quad \text{ (as $\hg\in A$, $\hh\in B$, and $A$ and $B$ $*$-commute)}\\
&=c\bullet(d\bullet 1 * \hh)*\hg.
\end{align*}

The common key $K$ is defined by
\[
K=d\bullet (c\bullet 1 * \hg)*\hh=c\bullet(d\bullet 1 * \hh)*\hg.
\]
Alice can compute the key $K$ using the right hand expression in the equation above; Bob can compute $K$ by computing the middle expression.

\section{The proposed attack}
\label{sec:attack}

Eve, the adversary, sees all public information, and also sees the elements $(p,g):=c\bullet 1 * \hg\in \Omega$ and $(q,h):=d\bullet 1*\hh\in\Omega$ that are transmitted between Alice and Bob. Eve's goal is to compute the shared key. 
Rather than attempting to compute Alice's private key material $c$ and $\hg$, or Bob's private key material $d$ and $\hg$, 
our attack will recover the shared key directly.

An overview of our attack is as follows. We first argue that the group $C$ can be replaced by a `linearised' version of $C$: this makes it easier to test membership in $C$. We then show that Eve does not need to compute Alice's or Bob's secret information in order to derive the shared key: more limited information suffices. (This information is specified in equations~\eqref{eqn:Eve_goal} and~\eqref{eqn:Eve_linear_combination} below.) Finally, we show how Eve can compute this information. 

For a group $H$ of $n\times n$ matrices over a field $\bF$, 
we write $\Alg(H)$ for the $\bF$-algebra generated by $H$~\cite{AlgSpan}. 
So $\Alg(H)$ is the set of all $\bF$-linear combinations of matrices in $H$. 
We write $\Alg^*(H)$ for the set of all invertible matrices in $\Alg(H)$. 

The groups $C$ proposed in the CBKAP satisfy $C=\Alg^*(C)$. More generally, we may assume that this is always the case.
To see this, first note $\Alg^*(C)$ and $D$ commute 
elementwise since every element of $\Alg^*(C)$ is a linear combination of elements in $C$.  Thus, $C$ may be replaced by $\Alg^*(C)$ to obtain a valid new instance of the protocol. Moreover, since $C\subseteq\Alg^*(C)$ the new instance of the protocol is more general than the original protocol: 
Alice can choose her matrix $c$ from the larger group $\Alg^*(C)$. So if we successfully recover the common key in 
\emph{every} new instance of the protocol, we can successfully recover the common key in the original instance. 

Thus, from now on, we assume that $C=\Alg^*(C)$. Let $\kappa_1,\kappa_2,\ldots,\kappa_r\in C$ be a basis for $\Alg(C)$. Such a basis is not difficult to compute, using standard techniques. Our assumption means that any invertible linear combination of the matrices $\kappa_i$ lies in $C$.

Let $P\trianglelefteq A$ be the \emph{pure subgroup of $A$}, defined by
\[
P=\{\,(\alpha,g)\in A: g=e\,\}.
\]
Then $\phi(P)$ is a subgroup of $\GL_n(\bF)$. Consider the subgroup $\Alg^*(\phi(P))$ of 
$\GL_n(\bF)$. 
Concretely, an element $\alpha'\in \Alg^*(\phi(P))$ is an invertible matrix of the form 
\[
\alpha'=\sum_{i=1}^k \ell_i\phi(\alpha_i)
\]
where $k\geq 0$, $\ell_i\in \bF$ and $(\alpha_i,e)\in P$.

Suppose that Eve finds elements $\tc\in C$, $\alpha'\in \Alg^*(\phi(P))$ and $(\ta,g)\in \hG$ such that
\begin{equation}
\label{eqn:Eve_goal}
(p,g)=\tc\bullet (\alpha',e) * (\ta,g).
\end{equation}
Moreover, suppose that Eve can find an elements $(\alpha_i,e)\in P$ 
and $\ell_i\in\bF$ such that
\begin{equation}
\label{eqn:Eve_linear_combination}
\sum_{i=1}^k\ell_i\phi(\alpha_i)=\alpha'.
\end{equation}
Then Eve can compute the common key, as follows. Firstly, she computes the matrix
\[
\beta'=\sum_{i=1}^k\ell_i\phi(^{h}\alpha_i).
\]
This computation is possible for Eve, since $h$ is part of the message $(q,h)=(d\phi(b),h)\in \Omega$ transmitted from Bob to Alice. 
Now, $(\alpha_i,e)\in P\leq A$, and so $(\alpha_i,e)$ $*$-commutes with all elements in $B$. Thus,
\[
(q\phi(^{h}\alpha_i),h)=d\bullet 1*(b,h)*(\alpha_i,e)=d\bullet 1*(\alpha_i,e)*(b,h).
\]
Eve then computes 
$\tc\bullet (q\beta',h)*(\ta,g)$. We claim that this is equal to the common key $K$. To see this, first note that 
\begin{align*}
(q\beta',h)&=\sum_{i=1}^k\ell_i(q\phi(^{h}\alpha_i),h)\\
&=\sum_{i=1}^k\ell_i(d\bullet 1*(\alpha_i,e)*(b,h))\\
&=\sum_{i=1}^k\ell_i(d\phi(\alpha_i)\phi(b),h)\\
&=(d\sum_{i=1}^k\ell_i\phi(\alpha_i)\phi(b),h)\\
&=(d\alpha'\phi(b),h)\\
&=d\bullet(\alpha',e)*(b,h).
\end{align*}
Hence
\begin{align*}
\tc\bullet (q\beta',h)*(\ta,g)&=\tc\bullet d\bullet(\alpha',e)*(b,h)*(\ta,g)\\
&=d\bullet \tc\bullet (\alpha',e)*(b,h)*(\ta,g)\\
&\quad\text{ (since $\tc\in C$ and $d$ centralises $C$)}\\
&=d\bullet \tc\bullet (\alpha',e)*(\ta,g)*(b,h)\\
&\quad\text{ (as $(\ta,g)\in A$ and $(b,h)\in B$ $*$-commute)}\\
&=d\bullet(p,g)*\hh\\
&=d\bullet (c\bullet 1 * \hg)*\hh\\
&=K.
\end{align*}
So it suffices to show that Eve can find elements $\alpha_i$,
$\ell_i$, $\ta$, $\tc$ and $\alpha'$ so that Equations~\eqref{eqn:Eve_goal} and~\eqref{eqn:Eve_linear_combination} are satisfied.

\paragraph{Precomputation stage: Find the $\alpha_i$.} Eve computes a collection of elements $(\alpha_i,e)$ such that the matrices $\phi(\alpha_i)$ form a basis of $\Alg(\phi(P))$. Once this is done, any 
$\alpha\in \Alg^*(\phi(P))$ can easily be written in the form~\eqref{eqn:Eve_linear_combination}. Eve does not need to know the messages $(p,g)$ and $(q,h)$ in this stage, so this stage can be carried out as a precomputation. Eve proceeds as follows.

Eve generates, as in~\cite{KalkaTeicherTsaban12},
short products $(a',g')$ of generators of $A$ such that the order $r$ of the permutation $g'$ is small ($n$ or less),
and computes $\alpha_1=(a',g')^r=(a'',e)$.
She repeats this procedure to generate $\alpha_2,\alpha_3\ldots$. (Eve may also take products of some of the previously generated elements $(\alpha_1,e)$, $(\alpha_2,e),\ldots ,(\alpha_{i-1},e)$ to define $(\alpha_i,e)$.) 
Eve stops when the dimension of the $\bF$-linear span of the matrices $\phi(\alpha_i)$ stops growing, and fixes a linearly independent subset of these matrices.

At the end of this process (relabelling after throwing linearly dependent elements $\phi(\alpha_i)$ away), Eve has $\alpha_1,\alpha_2,\ldots \alpha_r$ such that $\phi(\alpha_1),\phi(\alpha_2),\ldots,\phi(\alpha_r)$ are a basis for a subspace $V$ of $\Alg(\phi(P))$. Indeed, we expect (with high probability) that $V=\Alg(\phi(P))$. We assume that this is true from now on.

\paragraph{Stage~1: Find $\ta$.} Find a product of generators in $A$ whose second component is equal to $g$, using the method in~\cite{KalkaTeicherTsaban12}. Let $(\ta,g)$ be this product.
Define $\gamma\in \GL_n(\bF)$ by
\[
(\gamma,e)=(p,g)*(\ta,g)^{-1}.
\]

\paragraph{Stage~2: Find $\tc$.} Recall that Eve knows $\kappa_1,\kappa_2,\kappa_3,\ldots,\kappa_{r}\in C$ that form a basis of $\Alg(C)$. She finds (see below) field elements $x_1,x_2,\ldots,x_{r}\in\bF$ such that
\begin{align}
\label{eqn:V}
\quad&\gamma^{-1}(x_1\kappa_1+x_2\kappa_2+\cdots+x_{r}\kappa_{r})\in V,\text{ and}\\
\label{eqn:invertible}
\quad&x_1\kappa_1+x_2\kappa_2+\cdots+x_{r}\kappa_{r}\text{ is invertible}.
\end{align}
Set $\tc=x_1\kappa_1+x_2\kappa_2+\cdots+x_{r}\kappa_{r}$. Since $\tc$ is an invertible element of $\Alg(C)$, we see that $\tc\in C$.

To find a solution to Equations~\eqref{eqn:V} and~\eqref{eqn:invertible}, Eve randomly generates solutions $x_i$ that satisfy \eqref{eqn:V}, which is easy, as the conditions are linear. She stops when~\eqref{eqn:invertible} is also satisfied. We claim that the proportion of solutions to~\eqref{eqn:V} that satisfy~\eqref{eqn:invertible} is bounded below by $1-n/|\bF|$, which is a non-trivial proportion for the parameters that are proposed. The claim follows by applying the Invertibility Lemma~\cite[Lemma~9]{Tsaban12}, which states that the proportion of invertible matrices in any $\bF$-subspace of matrices over $\bF$ is at least $1-(n/|\bF|)$, provided that the subspace contains at least one invertible matrix. We note that the elements of the form $x_1\kappa_1+x_2\kappa_2+\cdots+x_{r}\kappa_{r}$ that satisfy~\eqref{eqn:V} are a subspace of matrices. So it remains to show that there exists an invertible element of this form. But let $x_1,x_2,\ldots,x_r\in \bF$ be such that $x_1\kappa_1+x_2\kappa_2+\cdots+x_{r}\kappa_{r}=c$. The elements $x_i$ exist since $c\in C\subseteq\Alg(C)$. Clearly $\tc=c$ is invertible. Moreover~\eqref{eqn:V} holds, because
we may show that $\gamma^{-1}c\in \phi(P)\subseteq V$ as follows. Firstly,
\[
(\gamma,e)=(p,g)*(\ta,g)^{-1}=(c\phi(a),g)*(\,^{g^{-1}}\!(\ta^{-1}),g^{-1})=(c\phi(a)\phi(\ta^{-1}),e),
\]
so $\gamma=c\phi(a)\phi(\ta^{-1})$ and therefore $\gamma^{-1}c=\phi(\ta)\phi(a)^{-1}$. And secondly, we see that $\phi(\ta)\phi(a)^{-1}=\phi(\ta a^{-1})\in \phi(P)$, since
\[
(\ta,g)(a,g)^{-1}=(\ta,g)(\,^{g^{-1}}\!a^{-1},g^{-1})=(\ta a^{-1},e)
\]
and $(\ta,g),(a,g)\in A$.

\paragraph{Stage 3: The remaining parameters.}

Eve sets $\alpha'=\tc^{-1}\gamma$. Since $(\alpha')^{-1}\in V$, we see that $\alpha'$ (being a power of $(\alpha')^{-1}$) also lies in $V$. So Eve can easily calculate coefficients $\ell_i$ such that
\[
\sum_{i=1}^k\ell_i\phi(\alpha_i)=\alpha'.
\]
Hence, Equation~\eqref{eqn:Eve_linear_combination} holds. We may also verify that Equation~\eqref{eqn:Eve_goal} holds:
\begin{align*}
\tc\bullet (\alpha',e)*(\ta,g)&=\tc\bullet(\tc^{-1}\gamma,e)*(\ta,g)=(\gamma,e)*(\ta,g)\\
&=((p,g)*(\ta,g)^{-1})*(\ta,g)=(p,g).
\end{align*}

\section{Experiments and conclusion}
\label{sec:implementation}

We have implemented our attack in Magma~\cite{BosmaCannon97}, running on one 2GHz core of a multi-core server. 
We used 5 sets of actual challenge parameters kindly provided by SecureRF. 
These parameters all used the values $|\bF|=256$ and $n=16$. The subgroup $A$ is specified by a generating set; each generator for $A$ is given as a word of length approximately 650 (notice the large parameter setting!) in the generating set $\mathcal{X}=\{\,(x_i(t),s_i):1\leq i\leq n-1\,\}$ for $\overline{M}\ltimes S_n$ defined in~\cite{AnshelAnshelGoldfeldLemieux06}. In all $5$ cases, our attack terminated successfully, producing the exact shared key. 
Our attack used less than 64MB of memory, and terminated in less than 8 hours. 
We would like to emphasise that our code is far from being optimised; we estimate an improvement in CPU time by a significant factor in an optimised version of the attack.

Let $B_n$ be a braid group on $n$ strands, and let $\sigma_1,\sigma_2,\ldots ,\sigma_{n-1}$ be the Artin generators for $B_n$. (See, for example,~\cite{KasselTuraev08} for an introduction to braid groups.) There is a homomorphism $\psi\colon B_n\rightarrow \overline{M}\ltimes S_n$ such that $\psi(\sigma_i)=(x_i,s_i)$ for $1\leq i\leq n-1$, which gives rise to the coloured Burau representation. 
Thus we could (and did) use standard routines for computing with braids in $B_n$, rather than dealing with words in $\mathcal{X}$ directly.

The most computationally intensive part of the attack is the computing of $\phi(a)$ where $(a,g)\in \overline{M}\ltimes S_n$ is given as a word in the generators of $A$. The long length of the generators in $A$ as words in $\mathcal{X}$ is the cause of difficulty here; we were computing with words of length approximately 20{,}000 in Stage~1 of our attack.

To decide when the precomputation stage should terminate, we use the criterion that the $\bF$-dimension of the algebra generated by the matrices $\phi(\alpha_i)$ should not grow when $4$ generators $(\alpha_i,e)$ in a row are considered.

Not surprisingly, this attack is highly parallelisable. We did not exploit this fact since for the actual parameters a single CPU core sufficed.

It remains open how to immunise the Algebraic Eraser against the presented cryptanalysis. The only hope seems to be to make the problem of expressing a permutation as a short product of given permutations difficult, by working with very carefully chosen distributions. However, for the intended applications, the computational constraints necessitate small values of $n$. In this case, Schreier--Sims methods solve this problem efficiently, no matter how the permutations are used. See the discussion around~\cite[Table~4]{KalkaTeicherTsaban12}.

\section*{Acknowledgement} The authors would like to thank Arkadius Kalka for providing code from the earlier attack~\cite{KalkaTeicherTsaban12} on the Algebraic Eraser, and for explaining how to use this code. The authors would also like to thank Martin Albrecht for various excellent editorial suggestions.

\end{document}